\documentclass{article}
\usepackage{amsmath}
\usepackage{amssymb}
\usepackage{graphicx}
\usepackage{bm}
\usepackage{enumerate}

\newtheorem{theorem}{Theorem}
\newtheorem{lemma}[theorem]{Lemma}
\newtheorem{observ}[theorem]{Observation}

\newtheorem{definition}[theorem]{Definition}
\def\qed{{\hfill$\Box$}}
\def\a{d^+}
\def\ad{\Delta^+}
\def\b{d^-}
\def\bd{\Delta^-}
\def\A{D^+}
\def\B{D^-}
\def\DD{\mathbf{(\A,\B)}}

\def\Proof{{\noindent \em Proof.}\ }
\newif\ifdeveloping

\ifdeveloping
\usepackage[notref,notcite]{showkeys}
\fi

\title{A simple Havel--Hakimi type algorithm to realize
graphical degree sequences of directed graphs\thanks{
PLE was partly supported by OTKA (Hungarian NSF), under contract Nos. AT048826 and K 68262.
IM was supported by a Bolyai postdoctoral stipend and OTKA (Hungarian NSF) grant F61730.
ZT was supported in part by the NSF BCS-0826958, HDTRA 201473-35045 and by Hungarian Bioinformatics MTKD-CT-2006-042794 Marie Curie Host Fellowships for Transfer of Knowledge. }
}

\author
    {P\'eter L. Erd\H{o}s\quad and \quad Istv\'an Mikl\'os \\
    {\small A. R\'enyi Institute of Mathematics, Hungarian
    Academy of}\\
    {\small Sciences, Budapest, PO Box 127, H-1364, Hungary}\\
    {{\small\tt \{elp\} \{miklosi\}@renyi.hu}}
\and
    \quad Zolt\'an Toroczkai \\
    {\small Interdisciplinary Center for Network Science
    and Applications } \\
    {\small and Department of Physics University of Notre Dame } \\
    {\small Notre Dame, IN, 46556, USA}\\
    {{\small\tt toro@nd.edu}}
    }

\date{\today}

\begin{document}

\maketitle

\begin{abstract}
\noindent One of the simplest ways to decide whether a given finite sequence of positive integers can arise as the degree sequence of a simple graph is the greedy algorithm of Havel and Hakimi. This note extends their approach to directed graphs. It also studies cases of some simple forbidden edge-sets. Finally, it proves a result which is useful to design an MCMC algorithm to find random realizations of prescribed directed degree sequences.
\end{abstract}

\section{Introduction}\label{sc:intro}

The systematic study of graphs (or more precisely the t{\em linear graphs}, as it was called in that time) began sometimes in the late forties, through seminal works by P. Erd\H{o}s, P. Tur\'an, W.T. Tutte, and others. One problem which received considerable attention was the existence of certain subgraphs of a given graph. For example such subgraph could be a perfect matching in a (not necessarily bipartite) graph, or a Hamiltonian cycle through all vertices, etc. Generally these substructures are called {\em factors}. The first couple of important results of this kind are due to W.T. Tutte who gave necessary and sufficient conditions for the existence of 1-factors and $f$-factors.

In the case of complete graphs, the existence problem of such factors is considerably easier. In particular, the existence problem of (sometimes simple) undirected graphs with given degree sequences  even admits simple greedy algorithms for its solution.

Subsequently, the theory was extended for factor problems of directed graphs as well, but the greedy type algorithm mentioned above, to the best knowledge of the authors, is missing even today.

In this paper we fill this gap: after giving a short and comprehensive (but definitely not exhausting) history of the $f$-factor problem (Section~\ref{sc:history}), we describe a greedy algorithm to decide the existence of a directed simple graph possessing the prescribed degree sequence (Section~\ref{sc:detail}). In Section~\ref{sc:Brual} we prove a consequence of the previous existence theorem, which is a necessary ingredient for the construction of edge-swap based Markov Chain Monte Carlo (MCMC) methods to sample directed graphs with prescribed degree sequence. Finally in Section~\ref{sec:constrained} we discuss a slightly harder existence problem of directed graphs with prescribed degree sequences where some vertex-pairs are excluded from the constructions. This result can help to efficiently generate all possible directed graphs with a given degree sequence.

\section{A brief history (of $f$-factors)}\label{sc:history}

For a given function $f:V(G)\rightarrow \mathbb{N}\cup \{0\}$, an {\em $f$-factor} of a given simple graph $G(V,E)$ is a subgraph $H$ such that $d_H(v)=f(v)$ for all $v\in V.$ One of the very first key results of modern graph theory is due to W.T. Tutte: in 1947 he gave a complete characterization of simple graphs with an
$f$-factor in case of $f\equiv 1$ (Tutte's 1-factor theorem, \cite{T47}). Tutte later solved the problem of the existence of $f$-factors for general $f$'s (Tutte's $f$-factor theorem, \cite{T52}). In 1954 he also found a beautiful graph transformation to handle $f$-factor problems via perfect matchings in bipartite graphs \cite{T54}. This also gave a clearly polynomial time algorithm for finding $f$-factors.

In cases where $G$ is a complete graph, the $f$-factor problem becomes easier: then we are simply interested in the existence of a graph with a given degree sequence (the exact definitions will come in Section~\ref{sc:detail}). In 1955 P. Havel developed a simple greedy algorithm to solve the {\em degree sequence problem} for simple undirected graphs (\cite{H55}). In 1960 P. Erd\H{o}s and T. Gallai studied the $f$-factor problem for the case of a complete graph $G$, and proved a simpler Tutte-type result for the degree sequence problem (see \cite{EG60}). As they already pointed out, the result can be derived directly form the original $f$-factor theorem, taking into consideration the special properties of the complete graph $G$, but their proof was independent of Tutte's proof and they referred to Havel's theorem.

In 1962 S.L. Hakimi studied the degree sequence problem in undirected graphs with multiple edges (\cite{H62}). He developed an Erd\H{o}s-Gallai type result for this much simpler case, and for the case of simple graphs he rediscovered the greedy algorithm of Havel. Since then this algorithm is referred to as the
{\em Havel--Hakimi algorithm}.

For directed graphs the analogous question of recognizability of a {\em bi-graphical-sequence} comes naturally. In this case we are given two $n$-element vectors $\mathbf {\a,\b}$ of non-negative integers. The problem is the existence of a directed graph on $n$ vertices, such that the first vector represents the out-degrees and the second one the in-degrees of the vertices in this graph. In 1957 D. Gale and H. J. Ryser independently solved this problem for simple directed graphs (there are no parallel edges, but loops are allowed), see \cite{G57, R57}. In 1958 C. Berge generalized these results for $p$-graphs where at most $p$ parallel edges are
allowed (\cite{B62}). (Berge calls the out-degree and in-degree together the {\em demi-degrees}.) Finally in 1973, the revised version of his book {\sl Graphs} (\cite{B73}) gives a solution for the $p$-graph problem, loops excluded. To show some of the afterlife of these results: D. West in his renowned recent textbook (\cite{W01}), discusses the case of simple directed graphs with loops allowed.

The analog of $f$-factor problems for directed graphs has a sparser history. {\O}ystein Ore started the systematic study of that question in 1956 (see \cite{O56,O56a}). His method is rather algebraic, and the finite and infinite cases - more or less - are discussed together. The first part developed the tools and proved the directly analog result of Tutte's $f$-factor problem for finite directed graphs (with loops), while the second part dealt with the infinite case.

In 1962 L.R. Ford and D.R. Fulkerson studied, generalized and solved the ``original" $f$-factor problem for a directed graph $\vec G$ (\cite{FF62}). Here lower and upper bounds were given for both demi-degrees of the desired subgraph (no parallel edges, no loops) with the original question naturally corresponding to equal lower and upper bounds. The solutions (as well as in Berge's cases) are based on network flow theory.

Finally, in a later paper Hakimi also proves results for bi-graphical sequences, however, without presenting a directed version of his original greedy algorithm (see \cite{H65}).

\section{A greedy algorithm to realize bi-graphical sequences}\label{sc:detail}

A sequence $\mathbf{d}=\{d_1,d_2,\ldots,d_n\}$ of nonnegative integers is called a {\em graphical sequence} if a simple graph
$G(V,E)$ exists on $n$ nodes, $V = \{v_1,v_2,\ldots,v_n \}$, whose degree sequence is $\mathbf{d}$. In this case we say that $G$ {\em realizes} the sequence $\mathbf{d}$. For simplicity of the
notation we will consider only sequences of strictly positive integers ($d_n > 0$) to avoid isolated points. The following, well-known result, was proved independently by V. Havel and S.L. Hakimi.
\begin{theorem}[Havel \cite{H55}, Hakimi \cite{H62}]\label{tm:HH}
There exists a simple graph with degree sequence $d_1>0,$
$d_2\ge \cdots \ge d_n > 0$ ($n \geq 3$) if and only if there exists one with degree sequence $d_2-1,\ldots,d_{d_1+1}-1,d_{d_1+2},\ldots,d_n$. (Note that there is no prescribed ordering relation between $d_1$ and the other degrees.)
\end{theorem}
This can be proved using a recursive procedure, which transforms any realization of the degree sequence into the form described in the Theorem \ref{tm:HH}, by a sequence of two-edge swaps.

\medskip\noindent A {\em bi-degree-sequence} (or {\em BDS} for
short) $\mathbf{(\a,\b)} = (\{\a_1,\a_2,\ldots,\a_n\}, \{\b_1, \b_2, $ $\ldots, \b_n\})$ of nonnegative integers is called a {\em bi-graphical sequence} if there exists a simple directed graph (digraph) $\vec G(V,\vec E)$ on $n$ nodes, $V = \{v_1,v_2, \ldots,v_n \}$, such that the out-degree and in-degree sequences together form $\mathbf{(\a,\b)}$. (That is the out-degree of vertex $v_j$ is $\a_j$ and its in-degree is $\b_j.$) In this case we say that $\vec G$ {\em realizes} our BDS. For simplicity, we will consider only sequences of strictly positive integer BDS's, that is each degree is $\ge 0$ and $\a_j+\b_j > 0$, to avoid isolated points.

Our goal is to prove a Havel--Hakimi type algorithm to realize bi-graphical sequences. To that end we  introduce the notion of normal order: we say that the BDS is in {\em normal order} if the entries satisfy the following properties: for each $i=1, \ldots ,n-2$ we either have $\b_i > \b_{i+1}$ or $\b_i=\b_{i+1}$ and $\a_i \ge \a_{i+1}$. Clearly, all BDS-s can be arranged into normal order. Note that we made no ordering assumption about node
$v_n$ (the pair $\a_n,\b_n$).

\begin{theorem}\label{tm:newHH}

Assume that the BDS $\mathbf{(\a,\b)}$  $($with $\a_j+\b_j > 0$, $j\in[1,n])$ is in normal order and $d_n^{+} > 0$ (recall: the out-degree of the last vertex is positive). Then $\mathbf{(\a,\b)}$ is bi-graphical if and only if the BDS
\begin{eqnarray}\label{eq:recurs}
&&\ad_k =\left\{\begin{array}{lll} \a_k & \quad \mbox{if }
       & k \neq n \\
       0 & \quad \mbox{if} & k = n\; ,
   \end{array}\right. \label{outred} \\
&&\bd_k =\left\{\begin{array}{lll} \b_k - 1 & \quad \mbox{if }
       & k \le \a_n\\
       \b_k & \quad \mbox{if } & k > \a_n\;\;\;,
   \end{array}\right. \label{inred}
\end{eqnarray}
with zero elements removed (those $j$ for which $\ad_j=\bd_j=0$) is bi-graphical.
\end{theorem}

\medskip\noindent Before starting the proof, we emphasize the similarity between this result and the original HH-algorithm. As in the undirected case, using Theorem~\ref{tm:newHH}, we can find in a greedy way a proper realization of graphical bi-degree sequences.

Indeed: choose any vertex $v_n$ with non-zero out-degree from the sequence, arrange the rest in normal order, then make $d_n^{-}$ connections from $v_n$ to nodes with largest in-degrees, thus constructing the out-neighborhood of $v_n$ in the (final) realization. Next, remove the vertices (if any) from the remaining sequence that have lost both their in- and out- degrees in the process, pick a node with non-zero out-degree, then arrange the rest in normal order. Applying Theorem 2 again, we find the final out-neighborhood of our second chosen vertex. Step by step we find this way the out-neighborhood of all vertices, while their in-neighborhoods get defined eventually (being exhausted by incoming edges). Note, that every vertex in this process is picked at most once, namely, when its out-neighborhood is determined by the Theorem, and never again after that.

\bigskip\noindent Our forthcoming proof is not the simplest, however, we use a more general setup to shorten the proofs of later results.

First, we define the partial order $\preceq$ among $k$-element vectors of increasing positive integers: we say $\mathbf{a} \preceq \mathbf{b}$ iff for each $j=1, \ldots ,k$ we have $a_j \le b_j.$

A {\em possible out-neighborhood} (or {\em PON} for short) of vertex $v_n$ is a $\a_n$-element subset of $V\setminus \{v_n\}$ which is a candidate for an out-neighborhood of $v_n$ in some graphical representation. (In essence, a PON can be any $\a_n$-element subset of $V\setminus \{v_n\}$ but later on we may consider some restrictions on it.) Let $A$ be a PON of $v_n.$
Then denote by $\bm{i}(A)$ the vector of the increasingly ordered subscripts of the elements of $A$. (For example, if $A=\{v_2,v_4,v_9\}$, then $\bm{i}(A) = (2,4,9)$.) Let $A$ and $B$ be two PONs of $v_n.$ We write:
\begin{equation}\label{def:left}
B \preceq A \quad \Leftrightarrow \quad \bm{i}_B \preceq \bm{i}_A\;.
\end{equation}
In this case we also say that $B$ is {\em to the left} of $A$. (For example, $B=\{v_1,v_2,v_6,v_7\}$ is to the left of $A=\{v_2,v_4,v_6,v_9\}$.)

\begin{definition} \label{reduced}
Consider a bi-graphical BDS sequence $(\mathbf{\a,\b})$ and let $A$ be a PON of $v_n$. The {\em $A$-reduced BDS} $\left (\mathbf{\a}\big|_A, \mathbf{\b}\big|_A \right)$ is defined as:
\begin{eqnarray}
    \a_k\big|_A & = &
    \left\{\begin{array}{lll} \a_k & \quad \mbox{if }
                & k \neq n \\
                0 & \quad \mbox{if} & k = n\; ,
            \end{array}\right. \label{4}\\
    \b_k\big|_A & = &
    \left\{\begin{array}{lll} \b_k - 1 & \quad
                \mbox{if } & k \in \bm{i}(A)\\
                \b_k & \quad \mbox{if } & k \not\in \bm{i}(A).
            \end{array}\right. \label{5}
\end{eqnarray}
\end{definition}
\noindent In other words, if $A$ is a PON in a BDS, then the reduced degree sequence $\left (\mathbf {\a}\big|_A, \mathbf{\b}\big|_A \right)$ is obtained by removing the out-edges of node $v_n$ (according to the possible out-neighborhood $A$). As usual, if for one subscript $k$ in the $A$-reduced BDS we have $\a_k\big|_A=\b_k\big|_A=0$ then the vertex with this index is to be removed from the bi-degree sequence.
\begin{lemma}\label{lm:left}
Let $(\mathbf {\a,\b})$ be a BDS, and let $A$ be a possible out-neighborhood of $v_n.$ Furthermore let $B$ be another PON with $B = A \setminus \{v_k\} \cup \{v_i\}$ where $\b_i \ge \b_k$ and in case of $\b_i = \b_k$ we have $\a_i \ge \a_k.$ Then if $\DD := \left (\mathbf {\a}\big|_A, \mathbf{\b}\big|_A \right)$ is bi-graphical, so is $\left (\mathbf {\a}\big|_B, \mathbf{\b}\big|_B \right)$.
\end{lemma}
\Proof Since our $A$-reduced BDS $\DD$ is bi-graphical, there exists a directed graph $\vec G$ which realizes the bi-degree sequence $\DD$. We are going to show that in this case there exists a directed graph $\vec G'$ which realizes the BDS $\left (\mathbf {\a}\big|_B, \mathbf{\b}\big|_B \right)$. In the following, $v_av_b$ will always mean a directed edge from node
$v_a$ to node $v_b$. Let us now construct the directed graph $\vec G_1$ by adding $v_n v$ directed edges for each $v\in A.$ (Since according to (\ref{4}), in $\DD$ the out-degree of $v_n$ is equal to zero, no parallel edges are created.) The bi-degree-sequence of $\vec G_1$ is $(\mathbf {\a,\b}).$ Our goal is to construct another realization $\vec G_1'$ of $(\mathbf {\a,\b})$ such that the deletion of the out-edges of $v_n$ in the latter produces the BDS $\left (\mathbf {\a}\big|_B, \mathbf{\b}\big|_B \right)$.

By definition we have $v_n v_k \in \vec E_1,$ (the edge set of $\vec G_1$) but $v_n v_i \not\in \vec E_1$. At first assume that there exists a vertex $v_\ell$ ($\ell \ne i,k,n$), such that $v_\ell v_i\in \vec E_1$ but $v_\ell v_k \not\in \vec E_1$.
(When $\b_i > \b_k$ then this happens automatically, however if
$\b_i=\b_k$ and $v_kv_i \in \vec E_1$ then it is possible that the in-neighborhood of $v_i$ and $v_k$ are the same - except of course $v_k$, $v_i$ themselves and $v_n$.) This means that now we can {\em swap} the edges $v_n v_k$ and $v_\ell v_i$ into $v_n v_i$ and $v_\ell v_k$. (Formally we create the new graph $\vec G_1'=(V,\vec E_1')$ such that $\vec E_1' = \vec E_1 \setminus \{v_n v_k, v_\ell v_i \} \cup \{v_n v_i, v_\ell v_k \} .$) This achieves our wanted realization.

\medskip\noindent Our second case is when $\b_i=\b_k,$ $v_kv_i \in \vec E_1$, and furthermore
\begin{equation}\label{eq:similar}
\mbox{for each }\ell \ne i,k,n \quad\mbox{we have}\quad
v_\ell v_i \in \vec E_1 \Leftrightarrow v_\ell v_k \in \vec E_1.
\end{equation}
It is important to observe that in this case $v_i v_k \not\in \vec E_1:$ otherwise some $v_\ell$ would not satisfy (\ref{eq:similar}) (in order to keep $\b_i=\b_k$).

Now, if there exists a subscript $m$ (different from $k, i, n$)
such that $v_i v_m \in \vec E_1$ but $v_k v_m \not\in \vec E_1,$
then we create the required new graph $\vec G_1'$ by applying the
following triple swap (or three-edge swap): we exchange the directed edges $v_n v_k, v_k v_i$ and $v_i v_m $ into $v_n v_i, v_i v_k$ and $v_k v_m$.

By our  assumption we have $\a_i \ge \a_k$. On one hand side if $\a_i > \a_k$ holds then due to the properties $v_kv_i \in \vec E$ and $v_i v_k \not\in \vec E,$ there exist at least two subscripts $m_1, m_2 \ne i,k$ such that $v_i v_{m_j} \in \vec E$ but $v_k v_{m_j} \not\in \vec E$ and at least one of them differs from $n$. Thus, when $\a_i > \a_k$, we do find such an $m$ for which the triple swap above can be performed.

The final case is when $\b_i = \b_k$ and $\a_i = \a_k$. If vertex $v_m$ does not exist, then we must have $v_i v_n \in \vec E_1$ (to keep $\a_i = \a_k$), and in this case clearly, $v_kv_n \notin \vec E_1$. Therefore, in this (final) case the graphical realization $\vec G_1$ has the properties $v_n v_k, v_kv_i, v_i v_n \in \vec E_1$ and $v_n v_i, v_i v_k, v_k v_n \not\in \vec E_1$.
Then the triple swap
\begin{equation}\label{eq:3swap}
\vec E_1' := \vec E_1\setminus \left \{ v_n v_k, v_kv_i, v_i v_n \right\}
\cup \left \{ v_n v_i, v_i v_k, v_k v_n \right \}
\end{equation}
will produce the required new graphical realization $\vec G_1'$. \qed

\begin{observ}\label{observ}
For later reference it is important to recognize that in all cases above, the transformations from one realization to the next one
happened with the use of two-edge or three-edge swaps.
\end{observ}

\begin{lemma}\label{lm:left2}
Let $(\mathbf {\a,\b})$ be an BDS and let $A$ and $C$ be two possible out-neighborhoods of $v_n.$ Furthermore assume that $C\preceq A$, that is $C$ is to the left of $A$. Finally assume that vertices in $A \cup C$ are in normal order. Then if $\left (\mathbf {\a}\big|_A, \mathbf{\b}\big|_A \right)$ is bi-graphical, so is $\left (\mathbf {\a}\big|_C, \mathbf{\b}\big|_C \right)$.
\end{lemma}
\Proof Since $C$ is to the left of $A$ therefore, there is a (unique) bijection $\phi: C \setminus A \rightarrow A \setminus C$ such that $\forall c \in C\setminus A$ : $\bm{i}(\{c\}) < \bm{i}(\{\phi(c)\})$ (the subscript of vertex $c$ is smaller than the subscript of vertex $\phi(c)$). (For example, if $A=\{v_4,v_5,v_6,v_7,v_8,v_9\}$ and $C=\{v_1,v_2,v_3,v_5,v_7,v_8\}$, then $C\setminus A = \{v_1, v_2,v_3\}$, $A\setminus C = \{v_4,v_6,v_9\}$, and $\phi$ is the map $\{v_1 \leftrightarrow v_4, v_2 \leftrightarrow v_6, v_3 \leftrightarrow v_9 \}$).

To prove Lemma~\ref{lm:left2} we apply Lemma~\ref{lm:left} recursively for each $c\in C \setminus A$ (in arbitrary order) to exchange $\phi(c)\in A$ with $c\in C$, preserving the graphical character at every step. After the last step we find that the sequence reduced by $C$ is graphical. \qed

\medskip\noindent{\bf Proof of Theorem~\ref{tm:newHH}}: We can easily achieve now the required graphical realization of $(\mathbf{\a,\b})$ if we use Lemma~\ref{lm:left2} with the current $A$, and $C=\{v_1, \ldots ,v_{\a_n}\}.$ We can do that since $(\mathbf{\a,\b})$ is in normal order, therefore the assumptions of Lemma~\ref{lm:left2} always hold.\qed

\section{A simple prerequisite for MCMC algorithms
to sample directed graphs with given BDS}\label{sc:Brual}

\medskip\noindent In practice it is often useful to choose uniformly a random element from a set of objects. A frequently used tool for that task is a well-chosen Markov-Chain Monte-Carlo method (MCMC for short). To that end, a graph is established on the objects and random walks are generated on it. The edges represent operations which can transfer one object to the other. If the Markov chain can step from an object $x$ to object $y$ with non-zero probability, then it must be able to jump to $x$
from $y$ with non-zero probability (reversibility). If the graph is connected, then applying the well-known Metropolis-Hastings algorithm, it will yield a random walk converging to the uniform
distribution starting from an arbitrary (even fixed) object.

To be able to apply this technique we have to define our graph (the Markov chain) $\mathcal{G}(\mathbf{\a}, \mathbf{\b})= (\mathcal{V},\mathcal{E})$. The vertices are the different possible realizations of the bi-graphical sequence $(\mathbf{\a}, \mathbf{\b}).$ An edge represents an operation consisting of a two or three-edge swap which transforms the first realization into the second one. (For simplicity, sometimes we just say {\em swap} for any of them.) We will show:
\begin{theorem}\label{tm:chain}
Let $\vec G_1, \vec G_2$ be two realizations of the same bi-graphical sequence $(\mathbf{\a}, \mathbf{\b}).$
Then there exists a sequence of swaps which transforms $\vec G_1$ into $\vec G_2$ through different realizations of the same bi-graphical sequence.

\end{theorem}

{\bf Remark}: In the case of undirected graphs the (original) analogous observation (needing only two-edges swaps) was proved by H.J. Ryser (\cite{R57}).

\medskip\noindent

\Proof We prove the following stronger statement:
\begin{itemize}\label{obs:e}
\item[(\maltese)] {\em there exists a sequence of at most $2e$  swaps which transform $\vec G_1$ into $\vec G_2$, where $e$ is the total number of out-edges in $(\mathbf{\a}, \mathbf{\b})$}
\end{itemize}
by induction on $e$. Assume that $(\maltese)$ holds for $e'<e$. We can assume that our bi-graphical sequence is in normal order on the first $n-1$ vertices and $\a _n>0.$ By Theorem \ref{tm:newHH}   there is a sequence $T_1$ ($T_2$) of  $d=\a_n$ many swaps which transforms $\vec G_1$  ($\vec G_2$) into  a $\vec G'_1$   ($\vec G'_2$) such that $\Gamma^+_{\vec G_1'}(v_n)=\{v_1,\dots, v_d\}$ ($\Gamma^+_{\vec G_2'} (v_n)=\{v_1,\dots, v_d\}$).

We consider now the directed graphs $\vec G_1''$ ($\vec G_2'$) derived from  directed graph $\vec G_1'$   (directed graph $\vec G_2'$) by deleting all out-neighbors of $v_n.$ Then both directed graphs realize the bi-graphical sequence $(\Delta^+, \Delta^-)$ which, in turn, satisfies relations (\ref{outred}) and (\ref{inred}). Therefore the total number of out-degrees is $e-d$ in both directed graphs, and by the inductive assumption there is a sequence $T$  of $2(e-d)$  many  swaps which transforms $\vec G''_1$ into $\vec G_2''$.

Now observe that if a swap transforms $\vec H$ into $\vec H'$, then the ``inverse swap" (choosing the same edges and non-edges and swap them) transforms $\vec H'$ into $\vec H$. So the swap sequence $T_2$ has an inverse $T_2'$ which transforms $\vec G_2'$ into $\vec G_2$. Hence  the sequence $T_1TT_2'$ is the required swap sequence: it transforms $\vec G_1$ into $\vec G_2$ and its length is at most $d+2(e-d)+d=2e$. \qed

\section{Is a BDS bi-graphical when one of its vertex's
out-neighborhood is constrained?}\label{sec:constrained}

In network modeling of complex systems (for a rather general reference see \cite{NBW06}) one usually defines a (di)graph with components of the system being represented by the nodes, and the
interactions (usually directed) amongst the components being
represented as the edges of this digraph. Typical cases include biological networks, such as the metabolic network, signal transduction networks, gene transcription networks, etc. The graph is usually inferred from empirical observations of the system
and it is uniquely determined if one can specify all the connections in the graph. Frequently, however, the data available from the system is incomplete, and one cannot uniquely determine this graph. In this case there will be a {\em set} ${\mathcal{D}}$ of (di)graphs satisfying the existing data, and one can be faced with:
\begin{enumerate}[{\rm (i)}]
\item finding a typical element of the class ${\mathcal{D}},$
\item or generating {\em all} elements of the class ${\mathcal{D}}$.
\end{enumerate}
(A more complete analysis of this phenomenon can be found in
\cite{K09}.) In Section~\ref{sc:Brual} we already touched upon
problem (i) when ${\mathcal{D}}$ is the class of all directed
graphs of a given BDS. The analogous Problem (ii) for undirected graphs was recently addressed in \cite{K09} which provides an
economical way of constructing all elements from ${\mathcal{D}}$.
In this Section we give a prescription based on the method from \cite{K09}, to solve (ii) for the case of all directed graphs with prescribed BDS. This is particularly useful from the point of view of studying the abundance of motifs in real-world networks: one needs to know first {\em all} the (small) subgraphs, or motifs, before we study their statistics from the data.

Before we give the details, it is perhaps worth making the following remark: Clearly, one way to solve problem (i) would be to first solve problem (ii), then choose uniformly from ${\mathcal{D}}$. However, in (those very small) cases when
reasonable answers can be expected for problem (ii), problem (i) is rather uninteresting. In general, however, (i) cannot be solved efficiently by the use of (ii).

\begin{figure}[htbp]
\bigskip \centering \vspace*{-0.5cm}
\includegraphics[width=2in]{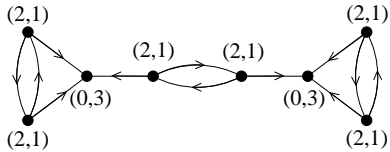}
\caption{\em This graph cannot be obtained by the directed Havel--Hakimi procedure. The integers indicate node
degrees.}\label{fig:nonghhD}
\end{figure}

We start the discussion of problem (ii) with pointing out that our new, directed Havel--Hakimi type algorithm is unable to generate all realization of a prescribed DBS (see Figure~\ref{fig:nonghhD}).

\noindent The situation is very similar to the non-directed case, see \cite{K09}. The directed HH-algorithm must start with a vertex with degree-pair $(2,1)$, therefore the two vertices of degree-pair $(0,3)$ must be out-neighbors of the same vertex - not for the graph in the Figure.

\medskip\noindent One possible way to overstep this shortage is to discover systematically all possible out-connections from a given vertex $v$ in all realizations of the prescribed graphical BDS.

We do not know a greedy algorithm to achieve this. The next best thing we can do is to develop a greedy algorithm to {\em decide} whether a given (sub)set of prescribed out-neighbors of $v$ would prevent to find a realization of the BDS containing those prescribed out-neighbors. In the following, we describe such a greedy algorithm. (It is perhaps interesting to note that this latter problem can be considered as a very special directed $f$-factor problem.)

To start, we consider a $\mathbf{(\a,\b)}$ bi-degree sequence together with a forbidden vertex set $F$ whose elements are not allowed to be out-neighbors of vertex $v_n.$ (Or, just oppositely, we can imagine that we already have decided that those vertices will become out-neighbors of $v_n$ and the BDS is already updated
accordingly. The forbidden vertex set governs only the out-neighbors, since in the process the in-neighbors are born ``automatically".) It is clear that $|F|+1 + \b_n \le n$ must hold for the existence of a graphical realization of this {\em $F$-restricted} BDS.

Assume that the vertices are enumerated in such a way that subset $F$ consists of vertices $v_{n-|F|}, \ldots ,v_{n-1}$ and vertices $V'=\{ v_1,  \ldots , v_{n-|F|-1}\}$ are in normal order. (We can also say that we apply a permutation on the subscripts accordingly.) Then we say that the BDS is in {\em $F$-normal order}.
\begin{definition} \label{reduced}
Consider a bi-graphical BDS sequence $(\mathbf{\a,\b})$ in $F$-normal order, and let $A$ be a PON. The {\em $A$-reduced BDS}
$\left (\mathbf{\a}\big|_A, \mathbf{\b}\big|_A \right)$ is defined as in (\ref{4}) and (\ref{5}), while keeping in mind the existence of an $F$ set to the right of $A$.
\end{definition}

\noindent In other words, if $A$ is a PON in an $F$-restricted BDS, then the reduced degree sequence $\left (\mathbf {\a}\big|_{A}, \mathbf{\b}\big|_{A} \right)$ is still obtained by removing the out-edges of node $v_n$ (according to the possible out-neighborhood $A$).

Finally, one more notation: let $(\mathbf{\a,\b})$ be a BDS, $F$ a forbidden vertex subset of $V$ and denote by $F[k]$ the set of the first $k$ vertices in the $F$-normal order.
\begin{theorem}\label{tm:kim}
Let $A$ be any PON in the $F$-restricted $(\mathbf{\a,\b})$ BDS,
which is in $F$-normal order. Then if the $A$-reduced BDS
$\left (\mathbf {\a}\big|_A, \mathbf{\b}\big|_A \right)$ is graphical, then the $F[\a_n]$-reduced BDS $\left (\mathbf {\a}\big|_{F[\a_n]}, \mathbf{\b}\big|_{F[\a_n]} \right)$ is graphical as well.
\end{theorem}
\Proof It is immediate: Lemma~\ref{lm:left2} applies. \qed

\medskip\noindent This statement gives us indeed a greedy way
to check whether there exists a graphical realization of the
$F$-restricted bi-degree sequence $(\mathbf{\a,\b})$: all we have
to do is to check only whether the $F[\a_n]$-reduced BDS $\left (\mathbf {\a}\big|_{F[\a_n]}, \mathbf{\b}\big|_{F[\a_n]} \right)$ is graphical.

Finally, we want to remark that, similarly to the indirected case, Theorem~\ref{tm:kim} is suitable to speed up the generation of all possible graphical realizations of a BDS. The details can be found in \cite{K09} which is a joint work of these authors with Hyunju Kim and L\'aszl\'o A. Sz\'ekely.

\subsection*{Acknowledgements}

The authors acknowledge useful discussions with G\'abor Tusn\'ady, \'Eva Czabarka and L\'aszl\'o A. Sz\'ekely and Hyunju Kim. ZT would also like to thank for the kind hospitality extended to him at the Alfr\'ed R\'enyi Institute of Mathematics, where this work was completed. Finally we want to express our gratitude to Antal Iv\'anyi for his editorial help.

\end{document}